\documentclass[11pt]{article} 

\usepackage[utf8]{inputenc}
\usepackage{amsmath}
\usepackage{amsfonts}
\usepackage{cases}
\usepackage{geometry} 
\usepackage{graphicx} 
\usepackage{subcaption}
\usepackage{color}

\numberwithin{equation}{section}

\newcommand{\nm}{\noalign{\smallskip}}
\newcommand{\beq}{\begin{equation}}
\newcommand{\eeq}{\end{equation}}
\newcommand{\eqnref}[1]{(\ref {#1})}
\newcommand{\ds}{\displaystyle}
\newcommand{\pf}{\noindent {\sl Proof}. \ }
\newcommand{\p}{\partial}
\newcommand{\pd}[2]{\frac {\p #1}{\p #2}}
\newcommand{\qed}{\hfill $\Box$}

 \def\p{\partial}
\def \Vh0{\stackrel{\circ}{V}_h} \def\to{\rightarrow}

\def\l{\label}  \def\f{\frac}  
\def\Box{\square}
\def\l|{\left|}
\def\r|{\right|}

\newtheorem{theorem}{Theorem}[section]

\newtheorem{lemma}[theorem]{Lemma}
\newtheorem{prop}[theorem]{Proposition}
\newtheorem{asump}[theorem]{Assumption}

\newtheorem{remark}{Remark}

\newcommand{\Bx}{\mathbf{x}}
\newcommand{\By}{\mathbf{y}}

\newcommand{\R}{\mathbb{R}}

\newcommand{\Scal}{\mathcal{S}}

\def\bx{{\mathbf{x}}}
\def\by{{\mathbf{y}}}
\def\bnu{{\mathbf{\nu}}}

\def\rd{{\mathrm{d}}}
\def\re{{\mathrm{e}}}
\def\ri{{\mathrm{i}}}

\def\cS{{\mathcal{S}}}

\def\be{\begin{equation}} \def\ee{\end{equation}}

\def\bea{\begin{eqnarray}}  \def\eea{\end{eqnarray}}

\def\beas{\begin{eqnarray*}} \def\eeas{\end{eqnarray*}}

\def\bn{\begin{enumerate}} \def\en{\end{enumerate}}

\def\bd{\begin{description}} \def\ed{\end{description}}

\title{Bloch waves in bubbly crystal near the first band gap: a high-frequency homogenization approach\thanks{\footnotesize The work of Hyundae Lee was supported by National Research Fund of Korea (NRF-2015R1D1A1A01059357, NRF-2017R1A4A1014735).  The work of Hai Zhang was partially supported by HK RGC GRF grant 16304517  and 
startup fund R9355 from HKUST.}}

\date{} 
\author{
Habib Ammari\thanks{\footnotesize Department of Mathematics, 
ETH Z\"urich, 
R\"amistrasse 101, CH-8092 Z\"urich, Switzerland (habib.ammari@math.ethz.ch).}  
\and Hyundae Lee\thanks{\footnotesize  Department of Mathematics, Inha University,  253 Yonghyun-dong Nam-gu,  Incheon 402-751,  Korea (hdlee@inha.ac.kr).}    \and Hai Zhang\thanks{\footnotesize 
Department of Mathematics, 
HKUST, Clear Water Bay, Kowloon, Hong Kong (haizhang@ust.hk).}}
   
\begin{document}
\maketitle
\begin{abstract}
This paper is concerned with the high-frequency homogenization of bubbly phononic crystals. It is a follow-up of the works [H. Ammari et al., Sub-wavelength phononic bandgap opening in bubbly media, J. Diff. Eq., 263 (2017), 5610--5629] which shows the existence of a sub-wavelength band gap. This phenomena can be explained by the periodic inference of cell resonance which is due to the high contrast in both the density and bulk modulus between the bubbles and the surrounding medium. In this paper, we prove that the first Bloch eigenvalue achieves its maximum at the corner of the Brillouin zone. Moreover, by computing the asymptotic of the Bloch eigenfunctions in the periodic structure near that critical frequency, we demonstrate that these eigenfunctions can be decomposed into two parts: one part is slowly varying and satisfies a homogenized equation, while the other is periodic across each elementary crystal cell and is varying. They rigorously justify, in the nondilute case, the observed super-focusing of acoustic waves in bubbly crystals near and below the maximum of the first Bloch eigenvalue and confirm the band gap opening near and above this critical frequency. 

\end{abstract}

\def\keywords2{\vspace{.5em}{\textbf{  Mathematics Subject Classification
(MSC2000).}~\,\relax}}
\def\endkeywords2{\par}
\keywords2{35R30, 35C20.}

\def\keywords{\vspace{.5em}{\textbf{ Keywords.}~\,\relax}}
\def\endkeywords{\par}
\keywords{bubble, sub-wavelength resonance, sub-wavelength phononic crystal,  Bloch theory, homogenization, metamaterial}

\section{Introduction}
This paper is devoted to the understanding of wave propagation in metamaterials (see for instance \cite{craster0}) which consists of sub-wavelength resonators arranged periodically in a background medium. These metamaterials differ from the usual photonic/phononic crystals (see for instance \cite{steve}) in the sense that their periods are much smaller than the free space wavelength of the functioning frequency of the materials. Note that in the later case, the periods are comparable to the wavelengths. Because of the sub-wavelength scale of the period, a homogenization theory is possible to describe the marcoscopic behavior of the materials, and this results in effective media having negative parameters such as negative mass, negative bulk modulus, negative electric permittivity, negative magnetic permeability, or double negative refractive index, or high contrast. The study of these metamaterials has drawn increasing interest nowadays because of their many important applications in fields such as superresolution, cloaking, and novel optic and phononic devices \cite{pendry2, pendry3}.  

There are many interesting mathematical works related to the homogenization theory for metamaterials, see for instance \cite{bouchitte, kohn, Lipton2, craster2}. Nevertheless, there is still much to understand about the wave propagation in these materials. Compared to the classic homogenization theory \cite{allaire, bensoussan, bloch,  milton}, the sub-wavelength  resonance \cite{phononic1, phononic2} 
produces strong interactions among the cells in the periodic structure and this can induce rich physics on the sub-wavelength scale which 
cannot be understood by the standard homogenization theory.  
Especially, one is interested in the high-frequency regime when the frequency is near the critical frequency where a sub-wavelength band gap of the periodic structure opens. Because this gap opens in the quasi-static regime (corresponding to a sub-wavelength scale), a homogenization theory is possible for such a structure. We note that the critical frequency usually occurs at the corner (or edge in the two dimensional case) of the Brillouin zone in three dimensions where one has typically anti-periodic Bloch eigenfunctions. In this frequency regime, the Bloch eigenfunctions vary on the microscale (the scale of the elementary crystal cell), and thus,  a homogenization theory which describes the macroscopic behaviour of the wave field seems impossible at first glance. This also makes the interpretation of the possible homogenization theory a perplexing task. On the other hand, the standard homogenization theory is applicable to the Bloch eigenfunctions which are near the center of the Brillouin zone and hence a much lower frequency regime. For this reason, the homogenization theory developed in this work is termed a high-frequency homogenization.

The bubbly media, because of the simplicity of the constituent resonant structure, the air bubbles,
become a natural model for such studies. It is known that a single bubble in the water possesses a sub-wavelength resonance which is called the Minnaert resonance \cite{Minnaert1933, H3a}. This resonance is due to due to the high contrast in both the density and bulk modulus between the bubbles and the surrounding medium. It makes the air bubble an ideal sub-wavelength resonator (the bubble can be two order of magnitude smaller than the wavelength at the resonant frequency). It is worth emphasizing that sub-wavelength resonators are the basic building block for metamaterials. As it can be seen from \cite{cooper2, cooper, zhikov}, the use of high contrast materials is not enough for the design of metamaterials, defined as a composite material with sub-wavelength focusing and sub-wavelength band-gap opening properties. 

We refer to \cite{leroy1, leroy2, leroy3, nature} for the experiments which motivated our series of studies of bubbles \cite{H3a,AFLYZ,Ammari_David,Ammari_Hai}. 
In \cite{Ammari_Hai}, using the fact that a single bubble can be well approximated by a monopole, and the point interaction approximation, we derived an effective medium theory for bubbly media consisting of dilute bubbles which may not be arranged periodically in a bounded domain. Our results show that, in the dilute case, near and below the Minnaert resonant frequency, the effective medium has high refractive index, which explains the super-focusing phenomenon observed in the experiment reported in \cite{leroy2}; while near and above the Minnaert resonant frequency, the effective medium is dissipative. 

Motivated by this work, we investigated the band structure of a bubbly phononic crystal which is made of periodically arranged bubbles in a homogeneous fluid \cite{AFLYZ}. We showed that there exists a sub-wavelength band gap in such a structure. This sub-wavelength band gap is mainly due to the cell resonance of the bubbles in the quasi-static regime and is quite different from the usual band gaps in photonic/phononic crystals where the gap opens at wavelength which is comparable to the period of the structure \cite{figotin, arma, Ammari2009_book}. We refer to \cite{Ammari_David} for the related work on bubbly metasurfaces which is a homogenization theory for a thin layer of periodically arranged bubbles mounted on a perfect reflection surface. In \cite{doublenegative}, a double-negative effective medium theory near the (hybridized) sub-wavelength resonances is obtained for an appropriate bubble volume fraction of randomly oriented bubble dimers. We also refer to \cite{caflish, caflish2} and the references therein for the other interesting related works on wave propagation in bubbly media.   

In this paper, based on the previous two works: effective medium theory for bubbly media in the dilute regime and the existence of a sub-wavelength band gap in the bubbly crystals, we further investigate the homogenization theory of the bubbly crystal near the frequency where the band gap opens. Our main approach is based on rigorous asymptotic analysis and layer potential techniques \cite{Ammari2009_book} which enables us to derive explicit formulas for the Bloch eigenfunctions. It is worth emphasizing that these formulas make both the homogenization theory and the justification of the super-focusing phenomenon in the nondilute case possible. 

We remark that our paper is related to the works \cite{craster1, craster2}. There are three major differences. (i) The homogenization in \cite{craster1, craster2} is concerned with the perturbation of the standing waves which are the Bloch eigenfunction at the edge of the Brillouin zone (in the two dimensional case). Our work is concerned with Bloch eigenfunctions near which the sub-wavelength band gap opens. We show that the band gap opens at the  corner (edge in two dimensions) of the Brillouin zone; 
(ii) The main approach in \cite{craster1, craster2} is based on a two-scale analysis, while our work relies on layer potential techniques; (iii)  The theory in \cite{craster1, craster2} may be restricted to the 
two dimensional case for certain structures, while our theory is applicable to any metamaterial where a sub-wavelength band gap exists. 
To sum up, this work complements the results of \cite{craster1, craster2}.

Finally, we would like to emphasize that the high frequency homogenization developed in the paper is very different from the classic one \cite{adams, allaire, allaire2, allaire3, arbogast, avila, bensoussan, bloch, babych, milton, kozlov, sparber,   briane, panasenko, rohan}. It has the following three new features: First, the cell structure can generate subwavelength resonance which has strong scattering effects near the resonant frequency. For the bubbly media considered here, the high contrasts in both the density and bulk modulus generate Minnaert resonance. Our homogenizatoin theory is developed around the frequencies such that Minnaert resonance can be excited. This is reflected in the scaling of the physical parameters adopted; Second, our homogenization exhibits phase transition phenomenon 
around the critical frequency where we derive the effective media theory. 
As is shown in the paper, below the critical frequency, the effective media is strongly dispersive, while above the critical frequency, the media becomes diffusive. This phase transition is mainly due to the anomalous scattering behaviour of resonators near their resonant frequencies. Third, our theory only hold for periodic media in the whole space where one can focus on the Bloch eigenfunction, but not for bounded domain with a large number of resonators. Indeed, because of the strong couple of resonators, a single resonator may dramatically change the scattered effect a system of coupled resonators near their resonator frequencies, see \cite{Ammari_Hai}. This also means that the uniform bound on the micro-field which are crucial for the success of the two-scale convergence theory does not hold. 
We believe that this homogenization theory provides the right mathematical framework for the investigation of wave propagation in metamaterials.

The paper is organized in the following way. In Section \ref{sec:setup}, we state the high frequency homogenization problem for the bubbly crystal, which we are interested in. Then in Section \ref{sec:pre}, we introduce some preliminaries on the layer potentials and quasi-periodic layer potential techniques. Next in Section \ref{sec:eigen}, we consider the normalized unit cell problem. We show that the first Bloch eigenvalue attains its maximum at the corner of the Brillouin zone. We also derive its asymptotic near that corner point. Finally in Section \ref{sec:homo}, we derive the high-frequency homogenization theory by analyzing the asymptotic of the Bloch eigenfunctions when the frequency is near the critical frequency which is the maximum of the first Bloch eigenvalue and where a band gap opens.

\section{Problem setup} \label{sec:setup}

We first describe the bubble phononic crystal under consideration. Let $Y$ be the unit cell $[-1/2,1/2]^3$ in $\mathbb{R}^3$, and let $D$ be a bounded and simply connected smooth domain contained in $Y$. 
The bubbles are periodically arranged with period $s>0$ in each direction. More precisely, let $D_s= sD$ be the domain occupied by the bubble in the unit cell $Y_s=[-\f{s}{2}, \f{s}{2}]^3$. Then the bubbles  
 occupy the domain $\cup_{n\in \mathbb{Z}^d} (s D+n)$. We denote by $\rho_b$ and $\kappa_b$ the density and the bulk modulus of the air inside the bubbles, respectively, and by $\rho$ and $\kappa$ the corresponding parameters for the background medium. Let $B_s= [-\f{\pi}{s}, \f{\pi}{s}]^3$
be the Brillouin zone corresponding to the periodic structure. 
The Bloch eigenvalues and eigenfunctions are solutions to the following $\alpha$-periodic equations in the cell $Y_s$ for each $\alpha \in B_s$:
\be \label{eq-scattering}
\left\{
\begin{array} {ll}
&\ds \nabla \cdot \f{1}{\rho} \nabla  u+ \frac{\omega^2}{\kappa} u  = 0 \quad \text{in} \quad Y_s \backslash D_s, \\
\nm
&\ds \nabla \cdot \f{1}{\rho_b} \nabla  u+ \frac{\omega^2}{\kappa_b} u  = 0 \quad \text{in} \quad D_s, \\
\nm
&\ds  u_{+} -u_{-}  =0   \quad \text{on} \quad \partial D_s, \\
\nm
& \ds  \f{1}{\rho} \f{\p u}{\p \bnu} \bigg|_{+} - \f{1}{\rho_b} \f{\p u}{\p \bnu} \bigg|_{-} =0 \quad \text{on} \quad \partial D_s,\\
\nm
&  e^{-i \alpha \cdot x} u  \,\,\,  \mbox{is periodic.}
  \end{array}
 \right.
\ee
Here, $\partial/\partial \bnu$ denotes the outward normal derivative and $|_\pm$ denote the limits from outside and inside $D$.  

Let
\begin{equation*} 
v = \sqrt{\frac{\kappa}{\rho}}, \quad v_b = \sqrt{\frac{\kappa_b}{\rho_b}}, \quad k= \frac{\omega}{v} \quad \text{and} \quad k_b= \frac{\omega}{v_b}
\end{equation*}
be respectively the speed of sound outside and inside the bubbles, and the wavenumber outside and inside the bubbles. We also introduce the dimensionless contrast parameter
\begin{equation*} 
\delta = \f{\rho_b}{\rho}. 
\end{equation*}
For bubbly media, we assume that $\delta \ll 1$, justifying the high contrast nature of the media.  In realistic setup, $\delta$ may be of the order of $10^{-3}$. On the other hand, we assume that 
$$
\f{k_b}{k}= \f{v}{v_b} =\sqrt{\f{\rho_b \kappa}{\rho \kappa_b}} = O(1),
$$
i.e., the wave numbers inside and outside the bubbles are comparable. 

It is known that \eqnref{eq-scattering} has nontrivial solutions for discrete values of $\omega$ which are called the Bloch eigenvalues. These eigenvalues can be arranged in the following increasing manner 
(see \cite{kuchment,reed}): 
$$  
0 \le \omega_{1,s}^\alpha \le \omega_{2,s}^\alpha \le \cdots.
$$
We denote the Bloch eigenfunction corresponding to the eigenvalue $\omega_{j,s}^\alpha$ by $u^{\alpha}_{j,s}$.

We have the following band structure of propagating frequencies for the given periodic structure:
$$ [0, \max_\alpha \omega_{1,s}^\alpha] \cup [ \min_\alpha \omega_{2,s}^\alpha, \max_\alpha \omega_{2,s}^\alpha] \cup  [ \min_\alpha \omega_{3,s}^\alpha, \max_\alpha \omega_{3,s}^\alpha] \cup \cdots. $$

In \cite{AFLYZ}, it is shown that there is a sub-wavelength band gap in the above band structure for fixed $s$ (say $s=1$) and sufficiently small $\delta$. More precisely, one has
$$
\omega_*^s:=\max_\alpha \omega_{1,s}^\alpha =O(\delta^{\f{1}{2}})<  \min_\alpha \omega_{2,s}^\alpha= O(1).
$$
In this paper, we investigate the asymptotic properties of the Bloch eigenfunctions in the bubbly crystal when the frequency $\omega$ is near the critical frequency $\omega_*^s$ where the sub-wavelength band gap opens for $s\ll 1$.  By a scaling argument, it can be shown that 
$$
\omega_*^s = \f{1}{s} \omega_*^1. 
$$ 
In order to fix the critical frequency in the limit when $s$ tends to zero, we rescale the contrast parameter $\delta$ as  follows:
\begin{equation} \label{scaling}
\delta = O(s^2). 
\end{equation}
Then the critical frequency remains of order one in the limiting process when $s$ tends to zero. Thus we are in the situation when the wavelength (of the free space) is of order one and the cell size is of order $s \ll 1$. As a result, one may develop a homogenization theory for such a structure. Indeed, in what follows, we shall show that the Bloch eigenfunctions $u^{\alpha}_{1,s}$, when the frequency is near the critical frequency $\omega_*^s$, can be decomposed into two parts: 
one part is slowly varying and satisfies a homogenized equation, while the other is periodic across each elementary crystal cell and is rapidly varying. This is very different from the usual homogenization where the second part is constant.
Our main approach is the layer potential techniques which we shall introduce in the next section. It is also worth emphasizing that the rapid variations of the second part of the solution justifies the super-focusing of acoustic waves in bubbly crystals experimentally observed in \cite{fink}.  Note also that the scaling in (\ref{scaling}) is well studied in the homogenization literature \cite{adams, arbogast, avila, babych, briane, panasenko, rohan}. Nevertheless, its importance here lies in the fact it induces a sub-wavelength resonance. As far as we know, the present paper is the first one to justify intriguing sub-wavelength properties of metamaterials at a critical sub-wavelength frequency corresponding to the maximum of the first Bloch eigenvalue.

\section{Preliminaries} \label{sec:pre}
We collect notations and some results regarding  the Green function and the quasi-periodic Green's function for the Helmholtz equation in three dimensions. We refer to \cite{Ammari2009_book} and the references therein for the details.

We  introduce the single layer potential $\mathcal{S}_{D}^{k} : L^2(\p D) \to H^1(\p D), H^1_{loc} (\R^3)$ associated with $D$ and the wavenumber $k$  defined by, $\forall \bx \in \R^3,$
$$
	  \mathcal{S}_{D}^{k} [\psi](\bx) :=  \int_{\p D} G^k(\bx, \by) \psi(\by) \rd \sigma(\by),
$$
where $$G^k(\bx, \by) :=  - \f{\re^{\ri k |\bx-\by|}}{4 \pi |\bx-\by|},$$ is the Green function of the Helmholtz equation in $\R^3$, subject to the Sommerfeld radiation condition. Here, $L^2(\p D)$ is the space of square integrable functions and $H^1(\p D)$ is the standard Sobolev space.

We also define the boundary integral operator $(\mathcal{K}_{D}^{k})^* : L^2(\p D) \to  L^2(\p D)$ by
$$
	(\mathcal{K}_{D}^{k})^* [\psi](\bx)  := \mbox{p.v.} \int_{\p D } \f{\p G_k(\bx, \by)}{ \p \nu(\bx)} \psi(\by) \rd \sigma(\by), \quad \forall \bx \in \p D.
$$
Here p.v. stands for the Cauchy principal value. 
 We use the notation $ \pd{}{\nu} \Big|_{\pm}$ indicating
$$ \pd{u}{\nu}\Big|_{\pm}(\bx)= \lim_{t \to 0^+} \langle \nabla u(\bx\pm t\nu(\bx)),\nu(\bx) \rangle,$$ 
with $\nu$ being the outward unit normal vector   to $\p D$. Then the following jump formula holds:
$$ \pd{}{\nu} \Big|_{\pm} \mathcal{S}_D^{k}[\phi](\bx) = \left( \pm \frac{1}{2} I + (\mathcal{K}_D^{k})^* \right)[\phi](\bx),\quad \mbox{a.e.}~\bx\in \p D.$$


We now define quasi-periodic layer potentials. 
Let $Y=Y_1$ be the unit cell $[-1/2,1/2]^3$. For $\alpha\in [-\pi,\pi]^3$, the function $G^{\alpha, k}$ is defined to satisfy
$$ (\triangle_\bx + k^2) G^{\alpha, k} (\bx,\by) = \sum_{n\in \mathbb{R}^3} \delta(\bx-\by-n) e^{\ri  n\cdot \alpha},$$
where $\delta$ is the Dirac delta function and  $G^{\alpha, k} $ is $\alpha$-quasi-periodic, {i.e.}, $e^{-\ri \alpha\cdot \bx} G^{\alpha, k}(\bx,\by)$ is periodic in $\bx$ with respect to $Y$. It is known that $G^{\alpha, k} $ can be written as
$$ G^{\alpha, k}(\bx,\by) = \sum_{n\in \mathbb{Z}^3} \frac{e^{\ri (2\pi n + \alpha)\cdot (\bx-\by)}}{k^2- |2\pi n + \alpha|^2},$$
if $k \ne |2\pi n + \alpha|$ for any $n\in \mathbb{Z}^3$. We remark that 
\begin{align} 
G^{\alpha,k}(\Bx,\By)=G^{\alpha,0}+ \sum_{\ell=1}^\infty k^{2\ell} G_\ell^{\alpha,\#}:= G^{\alpha,0}(\Bx,\By) - \sum_{\ell=1}^\infty k^{2\ell}\sum_{n\in\mathbb{Z}^d} \frac{e^{i(2\pi n + \alpha)\cdot (\Bx-\By)}}{|2\pi n +\alpha|^{2(\ell+1)}} \label{eq:defGk2}
\end{align}
when $\alpha \ne 0$, and $k \rightarrow 0$.

We are ready to define the quasi-periodic single layer potential $\mathcal{S}_D^{\alpha,k}$:
$$\mathcal{S}_D^{\alpha,k}[\phi](\bx) = \int_{\partial D} G^{\alpha,k} (\bx,\by) \phi(\by) d\sigma(\by),\quad \bx\in \mathbb{R}^3.$$
Then $\mathcal{S}^{\alpha,k}[\phi]$ is an $\alpha$-quasi-periodic function satisfying the Helmholtz equation $(\triangle + k^2)u=0$. In addition, one has the jump formula:
$$ 
\pd{}{\nu} \Big|_{\pm} \mathcal{S}_D^{\alpha,k}[\phi](\bx) = \left( \pm \frac{1}{2} I +( \mathcal{K}_D^{-\alpha,k} )^*\right)[\phi](\bx),\quad \mbox{a.e.}~\bx\in \p D,
$$
where $(\mathcal{K}_D^{-\alpha,k})^*$ is the operator given by
$$ 
(\mathcal{K}_D^{-\alpha, k} )^*[\phi](\bx)= \mbox{p.v.} \int_{\p D} \pd{}{\nu(x)} G^{\alpha,k}(\bx,\by) \phi(\by) d\sigma(\by).
$$
We remark that $\mathcal{S}_D^{0},~\mathcal{S}_D^{\alpha,0} : L^2(\p D) \rightarrow H^1(\p D)$ are invertible for $\alpha \ne 0$; see \cite{Ammari2009_book}. Moreover, the following decomposition holds for the layer potential $\mathcal{S}_D^{\alpha,k}$: 
\begin{equation} \label{series-s2}
\mathcal{S}_{D}^{\alpha,k} =  \cS_D^{\alpha, 0} + \sum_{\ell=1}^\infty k^{2\ell }\cS_{D,\ell}^{\alpha} \quad \text{with} \quad \cS_{D,\ell}^\alpha[\psi] := \int_{\p D} G_\ell^{\alpha,\#}(\bx - \by) \psi(\by) \rd \by,
\end{equation}
where the convergence holds in $\mathcal{B}(L^2(\p D), H^1(\p D))$, the set of linear bounded operators from $L^2(\p D)$ onto $H^1(\p D)$. 

Finally, we introduce the $\alpha$-quasi capacity of $D$, denoted by $\mbox{Cap}_{D,\alpha}$, 
$$
\mbox{Cap}_{D,\alpha}: = \int_{Y\setminus D} |\nabla u|^2, 
$$ 
where $u$ is the $\alpha$-periodic harmonic function in $Y\backslash \bar{D}$ with $u=1$ on $\p D$. For $\alpha\ne 0$,  we have $u(x) =\Scal_D^{\alpha,0} \left(\Scal_D^{\alpha,0}\right)^{-1}[1](x)$ and
$$
\mbox{Cap}_{D,\alpha} :=-\int_{\p D}  \left( \Scal_D^{\alpha,0}\right)^{-1}[1](y)d\sigma(y). 
$$ 
Moreover, we have a variational definition of $\mbox{Cap}_{D,\alpha}$. 
Indeed, let $C_{\alpha}^{\infty}(Y)$ be the set of $C^{\infty}$ functions in $Y$ which can be extended to $C^{\infty}$ $\alpha$-periodic functions in $\R^3$. Let $\mathcal{H}_{\alpha}$ be the closure of the set $C_{\alpha}^{\infty}(Y)$ in $H^1$, and let
$\mathcal{V}_\alpha:= \{  v\in \mathcal{H}_{\alpha} : v=1~\mbox{on }\p D\}$. Then we can show that
$$ 
\mbox{Cap}_{D,\alpha}= \min_{v\in\mathcal{V}_\alpha}\int_{Y\setminus D} |\nabla v|^2.
$$

\section{The Bloch eigenvalue in the unit period case} \label{sec:eigen}
We consider the bubbly phononic crystal when the period $s=1$ in this section.
For ease of notation, we write $\omega_{1,1}^\alpha= \omega_1^\alpha$ and $u_{1,1}^\alpha=u_1^\alpha$. 
We are interested in the point in the Brillouin zone where the maximum of $\omega_1^\alpha$ is achieved.  It is clear that due to time reversal symmetry, we have
$$
\omega_{1}^\alpha = \omega_{1}^{-\alpha}.
$$
On the other hand, 
$$
\omega_{1}^\alpha = \omega_{1}^{\alpha+(2\pi, 2\pi, 2\pi)}.
$$
Therefore,
$$
\omega_{1}^\alpha = \omega_{1}^{-\alpha+(2\pi, 2\pi, 2\pi)}.
$$
It then follows that $(\pi, \pi, \pi)$ is a critical point of $\omega_{1}^{\alpha}$. In what follows, we prove that $(\pi, \pi, \pi)$ is a maximum point with the following symmetry assumption on $D$.

\begin{asump} \label{asump0}
$D$ is symmetric with respect to planes $\{(x_1, x_2, x_3): x_j=0\}$, $j=1,2,3$.
\end{asump}

\begin{prop}
Suppose that Assumption \ref{asump0} holds. Let $\alpha^*:=(\pi,\pi,\pi)$. Then $\mbox{Cap}_{D,\alpha}$ and $\omega_1^\alpha$ attain their maxima at $\alpha=\alpha^*$.
\end{prop}
\pf
By the variational principle, we have the following characterization of $ \mbox{Cap}_{D,\alpha}$:
$$ 
\mbox{Cap}_{D,\alpha}= \min_{v\in\mathcal{V}_\alpha}\int_{Y\setminus D} |\nabla v|^2.
$$

Let $\mathcal{V}_{\alpha, 0}:= \{  v\in \mathcal{H}_{\alpha} : v=1~\mbox{on }\p D, ~~v=0\mbox{ on }\p Y\}$. Then $\mathcal{V}_{\alpha, 0} \subset \mathcal{V}_\alpha$ and

\begin{align*}
&\mbox{Cap}_{D,\alpha}= \min_{v\in{\mathcal{V}_\alpha}}\int_{Y\setminus D}  |\nabla v|^2 \le \min_{v\in{\mathcal{V}_{\alpha, 0}}}\int_{Y\setminus D}  |\nabla v|^2.
\end{align*}

Now, let $v_0$ be a  harmonic function satisfying 
\begin{equation}
\begin{cases} v_0 =1 \quad &\mbox{on} ~\p D,\\
v_0=0 &\mbox{on} ~\p Y.
\end{cases}
\end{equation}

Since $D$ is symmetric, $v_0$ is symmetric with respect to the planes $x_j=0$, $j=1,2,3$, and hence $v_0$ can be extended to an $\alpha^*$-quasi-periodic function (anti-periodic in each direction). 
By the characterization of $\mbox{Cap}_{D,\alpha^*}$, we have
$$
\int_{Y\setminus D}  |\nabla v_0|^2=\mbox{Cap}_{D,\alpha^*}.
$$
It follows that
\begin{align*}
&\mbox{Cap}_{D,\alpha} \le \min_{v\in\mathcal{Y}}\int_{Y\setminus D}  |\nabla v|^2 \le  \int_{Y\setminus D}  |\nabla v_0|^2=\mbox{Cap}_{D,\alpha^*}.
\end{align*}
Thus $ \mbox{Cap}_{D,\alpha}$ attains its maximum at $\alpha=\alpha^*$. 

The argument for $\omega_1^\alpha$ is similar to the one for $\mbox{Cap}_{D,\alpha}$.
By symmetry, we see that $u_1^{\alpha^*}$ is symmetric with respect to the planes $\{x_j=0\}$, $j=1,2,3$. Note that here I used the fact that the dimension of the Bloch eigenfunctions associated with $\omega_1^{\alpha^*}$ is one which is a consequence of the Gohberg-Sigal theory (see \cite{AFLYZ}). 
This combined with anti-periodicity yields that $u_1^{\alpha^*}$ is zero on $\p Y$. Let
 $$ 
\mathcal{H}_{\alpha, 0} :=\{  u \in \mathcal{H}_{\alpha}:  ~u=0\mbox{ on } \p Y\}.
$$
 
Then 
\begin{align*}
(\omega_1^\alpha)^2= \min_{u\in \mathcal{H}_\alpha} \frac{\int_Y \rho^{-1} |\nabla v|^2}{\int_Y \kappa^{-1} |v|^2}
 \le \min_{u\in\mathcal{H}_{\alpha, 0}} \frac{\int_Y \rho^{-1} |\nabla v|^2}{\int_Y \kappa^{-1} |v|^2} 
\le  \frac{\int_Y \rho^{-1} |\nabla u_1^{\alpha^*}|^2}{\int_Y \kappa^{-1} |u_1^{\alpha^*}|^2} = (\omega_1^{\alpha^*})^2.
\end{align*}
This completes the proof. \qed

\bigskip

\begin{remark}
We conjecture that the maximum of the first Bloch eigenvalue is achieved at the corner of the Brillouin zone for general periodic arrangements of sub-wavelength resonators under a similar symmetry assumption as Assumption \ref{asump0} on the cell structure. 
\end{remark}

In the sequel, we suppose that Assumption \ref{asump0} holds. To simplify the calculations and the presentation, we also assume 
\begin{asump} \label{asump1}
The wave speed inside the bubble is equal to the one outside, i.e., $v=v_b$.
\end{asump}

Let $u_1^\alpha$ be the $\alpha$-quasi-periodic propagating wave mode corresponding to $\omega_1^\alpha$, which is given in \cite{AFLYZ}  by 
\begin{equation}
u_1^\alpha=
\begin{cases}  \Scal_D^{\alpha,\omega_1^\alpha/ v} \left( \Scal_D^{\alpha,0}\right)^{-1} [1] +O(\delta^{1/2})\quad&\mbox{in}~ Y\setminus D,\\
\Scal_D^{\omega_1^\alpha/v_b} \left( \Scal_D^{0}\right)^{-1}[1] +O(\delta^{1/2})&\mbox{in}~D.
\end{cases}
\end{equation}

Since we have
\begin{align*}
\mathcal{S}_D^{\omega_1^\alpha/v_b}\left(\mathcal{S}_D^0\right)^{-1}[1]\approx\mathcal{S}_D^{0}\left(\mathcal{S}_D^0\right)^{-1}[1]
=\mathcal{S}_D^{\alpha,0}\left(\mathcal{S}_D^{\alpha,0}\right)^{-1}[1]
\approx\mathcal{S}_D^{\alpha,\omega_1^\alpha/v}\left(\mathcal{S}_D^{\alpha,0}\right)^{-1}[1]
\end{align*}
in $D$ up to a remainder of order of $O(\delta^{1/2})$, $u_1^\alpha$ can be approximated in $Y$ by
 \begin{equation}\label{u1_approx}
 u_1^\alpha=\Scal_D^{\alpha,\omega_1^\alpha/ v} \left( \Scal_D^{\alpha,0}\right)^{-1} [1] +O(\delta^{1/2}).
 \end{equation}

We shall investigate the behavior of  $\omega_1^\alpha$ near $\alpha=\alpha^*$.  
We first introduce some notation.
Let 
$$
G^{\tilde\alpha,1}\left(\Bx\right)= \sum_{n\in\mathbb{Z}^3} \frac{e^{i(2\pi n +\alpha^*)\cdot \Bx }}{k^2- |2\pi n + \alpha^*  |^2}\left( \frac{ |\tilde\alpha|^2}{k^2- |2\pi n + \alpha^*  |^2} +\frac{4 ((2\pi n+\alpha^*)\cdot\tilde\alpha)^2}{(k^2- |2\pi n + \alpha^*  |^2)^2} \right),
$$
and define the boundary integral operator
$$
\Scal_1^{\tilde\alpha}[\phi](\Bx):=\int_{\p D} G_1^{\tilde\alpha} (\Bx-\By) \phi(\By)d\sigma(\By).
$$

\begin{lemma}
The following holds
$$
\Scal_D^{\alpha^*+\epsilon\tilde\alpha,0}=
 e^{i \epsilon \tilde \alpha \cdot \Bx} \left( \Scal_{D}^{\alpha^*,0}+\epsilon^2 \Scal_1^{\tilde\alpha} +O(\epsilon^4) \right),
$$
where the $O(\epsilon^4)$ term is an operator from $L^2(\partial D)$ to $H^1(\partial D)$ whose operator norm is of order $\epsilon^4$.
\end{lemma}

\pf
Since
\begin{align*}&\frac{1}{k^2- |2\pi n + \alpha^* +\epsilon \tilde \alpha |^2}\\& =\frac{1}{k^2- |2\pi n + \alpha^*  |^2} \left( 1+\frac{\epsilon^2 |\tilde\alpha|^2}{k^2- |2\pi n + \alpha^*  |^2} +\frac{4\epsilon^2 ((2\pi n+\alpha^*)\cdot\tilde\alpha)^2}{(k^2- |2\pi n + \alpha^*  |^2)^2} \right) +O(\epsilon^3|\tilde\alpha|^3),
\end{align*}
using the symmetry in the summation, we can derive that
\begin{align*}
&G^{\alpha^*+\epsilon\tilde\alpha,k}\left(\Bx\right)=e^{i\epsilon\tilde \alpha\cdot \Bx}\sum_{n\in\mathbb{Z}^3} \frac{e^{i(2\pi n +\alpha^*)\cdot \Bx }}{k^2- |2\pi n + \alpha^* +\epsilon \tilde \alpha |^2} \nonumber\\
&= e^{i\epsilon\tilde \alpha\cdot \Bx}\sum_{n\in\mathbb{Z}^3} \frac{e^{i(2\pi n +\alpha^*)\cdot \Bx }}{k^2- |2\pi n + \alpha^*  |^2}\left( 1+\frac{\epsilon^2 |\tilde\alpha|^2}{k^2- |2\pi n + \alpha^*  |^2} +\frac{4\epsilon^2 ((2\pi n+\alpha^*)\cdot\tilde\alpha)^2}{(k^2- |2\pi n + \alpha^*  |^2)^2} \right) +O(\epsilon^4 |\tilde\alpha|^4)\nonumber \\
&= e^{i\epsilon\tilde \alpha\cdot \Bx} \left(G^{\alpha^*,k}\left(\Bx\right) +\epsilon^2 G_1^{\tilde\alpha}\left(\Bx\right) \right)+O(\epsilon^4|\tilde\alpha|^4). \label{gamma_exp}
\end{align*}
The lemma follows. \qed

\bigskip

Let $ \Lambda_D^{\tilde\alpha}$  be a quadratic function in $\tilde\alpha$  defined by
$$ 
\Lambda_D^{\tilde\alpha}:=\frac{1}{2} \int_{\partial D}  (\tilde \alpha \cdot \By)^2 \left( \Scal_D^{\alpha^*,0}\right)^{-1} [1]+ \left( \Scal_D^{\alpha^*,0}\right)^{-1}[(\tilde\alpha\cdot \Bx)^2]-2\left (\Scal_D^{\alpha^*,0}\right)^{-1}\Scal_1^{\tilde\alpha}\left(\Scal_D^{\alpha^*,0}\right)^{-1}[1] d\sigma(\By).
$$

\begin{lemma} \label{cap_approx} For every small $\epsilon>0$, it holds that
\begin{align*}
\mbox{Cap}_{D, \alpha^*+\epsilon\tilde\alpha}= \mbox{Cap}_{D, \alpha^*} +  \epsilon^2\Lambda_D^{\tilde\alpha} +O(\epsilon^4).
\end{align*}
Moreover, 
$\Lambda_D^{\tilde\alpha}$  is a negative semi-definite quadratic function of $\tilde \alpha$.  
\end{lemma}

\pf
Recall that 
$$
\mbox{Cap}_{D, \alpha^*+\epsilon\tilde\alpha}= -\int_{\p D}  \left( \Scal_D^{\alpha^*+\epsilon\tilde\alpha,0}\right)^{-1}[1](\By)d\sigma(\By).
$$

We solve
$
1= \Scal_D^{\alpha^*+\epsilon\tilde\alpha,0} [\phi]\left(\Bx\right)
$
for $\phi$. Since 
$\Scal_D^{\alpha^*+\epsilon\tilde\alpha,0}=
 e^{i \epsilon \tilde \alpha \cdot \Bx} \left( \Scal_{D}^{\alpha^*,0}+\epsilon^2 \Scal_1^{\tilde\alpha} +O(\epsilon^4) \right)$,
we get
\begin{align*}
\phi(\By)&= e^{i\tilde\alpha\cdot \epsilon \By} \left( \left( \Scal_D^{\alpha^*,0}\right)^{-1}-\epsilon^2\left (\Scal_D^{\alpha^*,0}\right)^{-1}\Scal_1^\alpha\left(\Scal_D^{\alpha^*,0}\right)^{-1} \right)[e^{-i\tilde\alpha \cdot \epsilon \Bx}](\By) +O(\epsilon^4)\\
&= e^{i\tilde\alpha\cdot \epsilon \By} \left[ \left( \Scal_D^{\alpha^*,0}\right)^{-1} [1]-\epsilon^2\left( \left( \Scal_D^{\alpha^*,0}\right)^{-1}[(\tilde\alpha\cdot \Bx)^2/2]-\left (\Scal_D^{\alpha^*,0}\right)^{-1}\Scal_1^{\tilde\alpha}\left(\Scal_D^{\alpha^*,0}\right)^{-1}[1] \right) \right](\By)\\
& +O(\epsilon^4).
\end{align*}

Using the expansion
$$
e^{i\tilde\alpha\cdot \epsilon \By} = 1 + i \epsilon \tilde\alpha\cdot \By - \frac{1}{2} \epsilon^2 (\tilde\alpha\cdot \By)^2 + \frac{-i \epsilon^3}{6}(\alpha\cdot  \By)^3 + O(\epsilon^4),
$$
and the symmetry in the integration, 
we have
\begin{align*}
\mbox{Cap}_{D, \alpha^*+\epsilon\tilde\alpha}&= - \int_{\partial D} \phi(\By) d\sigma(\By)\\
&= \mbox{Cap}_{D, \alpha^*} +  \epsilon^2\Lambda_D^{\tilde\alpha} +O(\epsilon^4).
\end{align*}
Finally, since $\mbox{Cap}_{D,\alpha}$ attains the maximum at $\alpha=\alpha^*$, $\Lambda_D^{\tilde\alpha}$  is a negative semi-definite quadratic function of $\tilde \alpha$.  This completes the proof. \qed

\bigskip

We now introduce a matrix $\{\lambda_{ij}\}$ such that
\begin{equation}
\frac{v_b^2}{|D|} \Lambda_D^{\tilde\alpha} := -\sum_{1\leq i, j \leq 3}\lambda_{ij}\tilde\alpha_i \tilde\alpha_j.
\end{equation}

It is clear that the matrix $\{\lambda_{ij}\}$ is symmetric and positive semi-definite.

\begin{lemma} \label{omega_exp} We have
\begin{align*}
(\omega_1^{\alpha^*+\epsilon \tilde \alpha}(\delta))^2 =\delta \omega_{M,\alpha^*}^2 -\epsilon^2 \delta \sum_{1\leq i, j \leq 3}\lambda_{ij} \tilde \alpha_j \tilde \alpha_j + \lambda_0(\alpha^*)\delta^2+O(\epsilon^4\delta+\epsilon^2\delta^2+\delta^3).
\end{align*}
\end{lemma}
\pf  
Recall that the following asymptotic formula of $\omega_1^\alpha$ in \cite{AFLYZ} holds:
\begin{equation}
(\omega_1^\alpha(\delta))^2= \frac{\delta v_b^2 \mbox{Cap}_{D,\alpha}}{|D|}+O(\delta^2):=\delta(\omega_{M,\alpha})^2+O(\delta^2),
\end{equation}
where $\omega_{M,\alpha} =  = \sqrt{\f{v_b^2\mbox{Cap}_{D,\alpha} }{|D|}}$, and the $O(\delta^2)$ term is uniformly of order $\delta^2$ with respect to $\alpha$ away from $0$. A further asymptotic expansion implies that the $O(\delta^2)$ term can be further decomposed as 
$$
O(\delta^2)= \lambda_0(\alpha)\delta^2 + O(\delta^3),
$$
where $O(\delta^3)$ denotes some uniformly bounded remaining term. Since $\alpha^*$ is a critical point of $\omega_1^\alpha$ regardless of $\delta$, we have $\lambda_0(\alpha)=\lambda_0(\alpha^*)+O(|\alpha-\alpha^*|^2)$.

Finally, we obtain that
\begin{align*}
(\omega_{M, \, \alpha^*+\epsilon\tilde\alpha})^2= \frac{ v_b^2 \mbox{Cap}_{D, \alpha^*+\epsilon\tilde\alpha}}{|D|}
&= (\omega_{M,\alpha^*})^2  + \frac{\epsilon^2 v_b^2}{|D|}\Lambda_D^{\tilde\alpha} +O(\epsilon^4).
\end{align*}
This completes the proof. \qed

\section{Homogenization of the Bloch eigenfunction} \label{sec:homo}
We now consider the asymptotic behaviours of the Bloch eigenfunction in the bubbly crystal with period $s$.
We first present the relation between the Bloch eigenvalues and eigenfunctions at different scales, which follows from a scaling argument.

\begin{lemma} \label{eq-scaling}
Let $\omega_1^{\alpha}$ and $u_1^{\alpha}$ be the Bloch eigenvalue and eigenfunction of the bubbly crystal with period $1$, then the bubbly crystal with period $s$ has Bloch eigenvalue 
$$\omega_{1,s}^{\alpha/s}=\f{1}{s}\omega_{1}^{\alpha}$$
and eigenfunction
$$
u_{1,s}^{\alpha/s}(\Bx) = u_1^{\alpha}\left (\f{\Bx}{s}\right).
$$
\end{lemma}

As a consequence, we see that
$$
 \omega_*^s= \max_{\alpha} \omega_{1,s}^{\alpha} = \frac{1}{s}\max_{\alpha} \omega_{1}^{\alpha}= \frac{1}{s} \omega_1^{\alpha^*}.
$$
Moreover, the maximum is attained at the corner point $\frac{1}{s}\alpha^*$ of the Brillouin zone $B_s$. 
We are interested in the behaviour of the Bloch eigenfunction $u_{1,s}^{\alpha}$ for $\alpha$ near $\alpha^*/s$.

\begin{lemma}
\begin{equation}
u_{1,s}^{\alpha^*/s + \tilde\alpha}(\Bx)= e^{i\tilde\alpha \cdot\Bx} S \left(\frac{\Bx}{s}\right) + O(s^2+\delta^{1/2}),
\end{equation}
where the function $S$ is defined by 
\begin{equation} \label{defS}
S(\Bx)= \Scal_{D}^{\alpha^*, 0}  \left[\left( \Scal_D^{\alpha^*,0}\right)^{-1} [1]\right] (\Bx).
\end{equation}
\end{lemma}

\pf
Note that $u_{1,s}^{\alpha^*/s + \tilde\alpha}(\Bx) =u_1^{\alpha^*+s\tilde\alpha}(\Bx/s)$. By  \eqnref{u1_approx}, we have 
\[
u_1^{\alpha^*+s\tilde\alpha}(\Bx/s) = \Scal_D^{\alpha^*+s \tilde \alpha, v^{-1}\omega_1^{\alpha^*+s\tilde\alpha}}\left [ \left(\Scal_D^{\alpha^*+s \tilde \alpha,0}\right)^{-1}[1]\right] \left( \frac{\Bx}{s}\right)+O(\delta^{1/2}).
\]
Since $\omega_1^{\alpha^*+s\tilde\alpha} = O(\delta^{1/2})$, using (\ref{series-s2}), we further obtain
\begin{equation}\label{asymp_10}
u_1^{\alpha^*+s\tilde\alpha}(\Bx/s)=\Scal_D^{\alpha^*+s \tilde \alpha, 0}\left [ \left(\Scal_D^{\alpha^*+s\tilde\alpha,0}\right)^{-1}[1]\right] \left( \frac{\Bx}{s}\right)+O(\delta^{1/2}).
\end{equation}
On the other hand, in the proof of Lemma \ref{cap_approx}, we showed that
\begin{equation}
\left( \Scal_D^{\alpha^*+s\tilde\alpha,0}\right)^{-1} [1]\left(\frac{\By}{ s}\right)=e^{i\tilde\alpha\cdot \By} \left( \Scal_D^{\alpha^*,0}\right)^{-1} [1]\left(\frac{\By}{ s}\right)+O(s^2).\label{asymp_11}
\end{equation}

It follows from \eqnref{asymp_10} and \eqnref{asymp_11} that
\begin{align*}
u_1^{\alpha^*+s\tilde\alpha}(\Bx/s) &= \int_{s \p D} e^{i \tilde \alpha \cdot (\Bx-\By) } G^{\alpha^*, 0} \left(\frac{\Bx-\By}{s}\right) e^{i\tilde\alpha\cdot  \By} \left( \left( \Scal_D^{\alpha^*,0}\right)^{-1} [1]\right)\left(\frac{\By}{s}\right) s^{-2} d\sigma(\By)\\
& +O(s^2+\delta^{1/2})\\
&= e^{i\tilde\alpha \cdot\Bx} \Scal_{D}^{\alpha^*, 0}  \left[\left( \Scal_D^{\alpha^*,0}\right)^{-1} [1]\right] \left(\frac{\Bx}{s}\right)+O(s^2+\delta^{1/2}).
\end{align*}

Then the lemma follows.
\qed

\bigskip

We note that the function $S$ defined by (\ref{defS}) is a piecewise harmonic function with $S=1$ on $\p D$ and $S=0$ on $\p Y$. $S \left(\frac{\Bx}{s}\right)$ varies on the small scale $s$ and describes the microscopic behaviour of the $u_{1,s}^{\alpha^*/s + \tilde\alpha}$. On the other hand, 
the function $e^{i\tilde\alpha\cdot \Bx}$ represents the macroscopic behaviour of $u_{1,s}^{\alpha^*/s + \tilde\alpha}$. We now derive a homogenized equation for this macroscopic field near the critical frequency $\omega_*^s$.

We first recall from \cite{AFLYZ} that 
$$
\omega_{1}^{\alpha^*} = C \sqrt{\delta} + O(\delta^{3/2})
$$
for some positive constant $C$ of order one which depends on $D$, $v_b$ and $\alpha^*$. In order to keep the critical frequency $\omega_*^s$ on a fixed order when the cell size $s$ tends to zero, we assume that

\begin{asump} \label{asump2}
$\delta = \mu s^2$ for some positive constant $\mu$.
\end{asump}

As a result, we have
$$
\max_{\alpha} \omega_{1,s}^{\alpha} = \omega_*^s = \frac{1}{s}\omega_{1}^{\alpha^*} =O(1).
$$

Now, suppose that  $\omega$ is near $\omega_*^s$. We need to find the corresponding Bloch eigenfunctions, or $\tilde \alpha$ so that 
$$
\omega^2 = \omega_{1,s}^{\alpha^*/s + \tilde\alpha}. 
$$

Using Lemmas \ref{omega_exp} and \ref{eq-scaling},  we have
\begin{equation} \omega_*^2-\omega^2 = \delta\sum_{1\leq i, j \leq 3}\lambda_{ij} \tilde \alpha_i\tilde\alpha_j+O(s^4),\label{omega_eq}\end{equation}
where $\omega_*:= \omega_*^s$. 

Let $$\omega_*^2 -\omega^2 = \beta \delta$$
for some constant $\beta$. Then we have
$$
\sum_{1\leq i, j \leq 3}\lambda_{ij} \tilde \alpha_i\tilde\alpha_j =\beta + O(s^2),
$$
which shows that all the solutions $\tilde \alpha$ lie approximately on the ellipsoid defined by 
$$
\{\alpha=(\alpha_1, \alpha_2, \alpha_3):  \sum_{1\leq i, j \leq 3}\lambda_{ij} \alpha_i \alpha_j =\beta\}.
$$ 
It also implies that the plane wave  $u(\Bx):=e^{i\tilde\alpha\cdot \Bx}$ satisfies to a leading order term
\begin{equation}
\sum_{1\leq i, j \leq 3}\lambda_{ij} \partial_i \partial_j \hat u (\Bx)+ \beta \hat u(\Bx)= 0.
\end{equation} 

We now consider two cases.
Case I: $\beta>0$ and Case II: $\beta <0$.
 
In the first case, $\tilde{\alpha}$ is a well-defined vector in $\R^3$, which justifies the existence of true Bloch eigenfunctions. However, in the second case $\tilde{\alpha}$ is a pure imaginary vector, and hence is not a Bloch eigenfunction. The associated function $e^{i\tilde\alpha\cdot \Bx}$ either grows or decays exponentially along certain direction. This justifies that a band gap opens at the critical frequency $\omega_*$.

To conclude, we have the following main result on the homogenization theory for the bubbly crystals. 

\begin{theorem}
Under Assumptions \ref{asump0}, \ref{asump1} and \ref{asump2}, for frequencies in a small neighborhood  of maximum of the first Bloch eigenvalue, say, $\omega_*^2 -\omega^2 = O(s^2)$, the following asymptotic of Bloch eigenfunction $u_{1,s}^{\alpha^*/s + \tilde\alpha}$ holds:
$$
u_{1,s}^{\alpha^*/s + \tilde\alpha}(\Bx) = e^{i \tilde\alpha \cdot \Bx} S \left(\frac{\Bx}{s}\right) + O(s),
$$
where the macroscopic field $e^{i\tilde\alpha \cdot\Bx}$ satisfies the following equation 
\begin{equation}
 \sum_{1\leq i, j \leq 3}\lambda_{ij} \partial_i \partial_j \hat u (\Bx)+ \frac{\omega_*^2 -\omega^2}{\delta}\hat u(\Bx)= 0.
\end{equation}
which can be viewed as the homogenized equation for the 
bubbly phononic crystal, while the microscopic field is periodic and varies on the scale of $s$. 
\end{theorem}

It is clear that the homogenized medium is very dispersive below and near the critical frequency $\omega_*$. Moreover, the microscopic oscillations of the field at the period of the crystal justify the super-focusing phenomenon \cite{fink}. This mechanism differs from the one based on the high contrast effective medium theory derived in \cite{Ammari_Hai}, which is valid only in the dilute case.  
 

%
%

\end{document}